\documentclass[11pt]{article}
\usepackage[centertags]{amsmath}
\usepackage{amsmath}
\usepackage{amssymb}
\usepackage{amsthm}

\numberwithin{equation}{section}
\newcommand{\dint}{\displaystyle\int}
\newcommand{\1}{1\!\!1}
\newtheorem{definition}{Definition}[section]
\newtheorem{theorem}{Theorem}[section]
\newtheorem{proposition}{Proposition}[section]
\newtheorem{lemma}{Lemma}[section]
\newtheorem{remark}{Remark}[section]
\newtheorem{corollary}{Corollary}[section]

\newcommand{\ind}{\mathbf{1}}

\newcommand{\eps}{\varepsilon}

\def \R{\mathbb{R}}

\def \N{\mathbb{N}}
\def \E{\mathbb{E}}

\def \L{\mathbb{L}}

\def \bf{\textbf}
\def \it{\textit}

\def \bop {\noindent\textbf{Proof }}
\def \eop {\hbox{}\nobreak\hfill
\vrule width 2mm height 2mm depth 0mm
\par \goodbreak \smallskip}

\def \F{\mathcal{F}}

\def \bop {\noindent\textbf{Proof }}
\def \eop {\hbox{}\nobreak\hfill
\vrule width 2mm height 2mm depth 0mm \par \goodbreak \smallskip}

\usepackage{amssymb}
\usepackage{amsthm}
\usepackage{amsmath}
\usepackage{latexsym}
\usepackage[T1]{fontenc}
\usepackage[latin1]{inputenc}

\def \R{\mathbb{R}}

\def \N{\mathbb{N}}
\def \E{\mathbb{E}}

\def \L{\mathbb{L}}

\def \bf{\textbf}
\def \it{\textit}

\def \bop {\noindent\textbf{Proof.}}
\def \eop {\hbox{}\nobreak\hfill
\vrule width 2mm height 2mm depth 0mm
\par \goodbreak \smallskip}

\def \F{\mathcal{F}}

\def \eop {\hbox{}\nobreak\hfill \vrule width 2.0mm height 1.8mm depth 0mm
\par \goodbreak \smallskip}

\numberwithin{equation}{section}

\def \R{\mathbb{R}}

\def \N{\mathbb{N}}
\def \E{\mathbb{E}}

\def \L{\mathbb{L}}

\def \bf{\textbf}
\def \it{\textit}

\def \bop {\noindent\textbf{Proof }}
\def \eop {\hbox{}\nobreak\hfill
\vrule width 2mm height 2mm depth 0mm
\par \goodbreak \smallskip}

\def \F{\mathcal{F}}

\begin{document}

%%%%%%%%%%%%%%%%%%%%%%%%%%%%%%%%%%%%%%%%%%%%%%%%%%%%%%%%%%%%%%%

\title{ Stochastic Optimal Control and BSDEs with\\ Logarithmic Growth\thanks{The research leading to
these results has received funding from the European Community's FP
7 Programme under contract agreement PITN-GA-2008-213841, Marie
Curie ITN "Controlled Systems".} }
\author{KHALED BAHLALI \thanks{%
IMATH, UFR Sciences, USVT, B.P. 132, 83957 La Garde Cedex, France.
e-mail: bahlali@univ-tln.fr } and BRAHIM EL ASRI \thanks{%
Institut f\"{u}r Stochastik Friedrich-Schiller-Universit\"{a}t Jena
Ernst-Abbe-Platz 2 07743 Jena, Germany; e-mail:
brahim.el-asri@uni-jena.de }}
\date{}
\maketitle \noindent {\bf{Abstract}} In this paper, we study the
existence of an optimal strategy for the stochastic control of
diffusion in general case and a saddle-point for zero-sum stochastic
differential games. The problem is formulated as an extended BSDE
with logarithmic growth in the $z$-variable and terminal value in
some $L^p$ space. We also show the existence and uniqueness of
solution of this
BSDE.\\
%\end{abstract}
\\
\noindent \text{AMS 2000 Classification subjects:} 60G40, 62P20,
91B99, 91B28, 35B37, 49L25.

\noindent {${Keywords}:$} Backward stochastic differential
equations, Stochastic control, Zero-sum stochastic differential
games.

\medskip

\section{Introduction}
\noindent In this paper we study BSDE with the applications to
stochastic
control and stochastic zero-sum differential games.\\
 \indent We consider a backward stochastic differential
equation (BSDE) with generator $\varphi$ and terminal condition
$\xi$

\begin{equation}\label{zlogz1} \displaystyle Y_t = \xi + \int_t^T
\varphi(s,Y_s,Z_s) ds - \int_t^T Z_s dB_s,\qquad
t\in[0,T]\end{equation} where $(B_t)_{t\geq0}$ is a standard
Brownian motion. Such equations have been extensively studied since
the first paper of E. Pardoux and S. Peng \cite{PP}. We will
consider the case when $\varphi$ is allowed to have logarithmic
growth $(\vert z\vert\ln^{\frac{1}{2}}(|z|))$ in the z-variable.
Moreover, we will
allow $\xi$ to be unbounded.\\
\bf{1)}\cite{HLP95} showed the existence of an optimal stochastic
control in the stochastic control of diffusions, in the case where
the drift term of equation $f$ which defines the controlled system
is bounded. In the same bounded case the existence of a saddle-point
for a zero-sum stochastic differential game can be proved in a
similar way.\\
\bf{2)}\cite{HLP295} established the existence of an optimal
stochastic control in the stochastic control of diffusions, in the
case where
 the running reward function $h$ is bounded. In the
same bounded case the existence of a saddle-point for a zero-sum
stochastic differential game can be proved in a similar way.\\

\indent Our aim in this work is to relax the boundedness assumption
on drift term of equation $f$ functionals and the running reward
function $h$. Therefore the main objective of our work, and this is
the novelty of the paper, is to show the existence of an optimal
strategy for the stochastic control of diffusion. The main idea
consists to showed the existence and uniqueness of the solution of
BSDE \ref{zlogz1} and characterize the value function as a solution
of BSDE.\\

This paper is organized as follows: In Section 2, we present the
assumptions and we formulate the problem. In Section 3, we give the
the main result on existence and uniqueness of the solution of BSDE
\ref{zlogz1}. In Section 4, we state some estimates of the solutions
from which we derive some integrability properties of the solution.
In Section 5, we give  estimate between two solutions and the proof
of Theorem \ref{unique} and \ref{stability}. In Section 6, we
introduce the optimal stochastic control problem and we give the
connection between optimal stochastic control problem and the
zero-sum stochastic differential games and the BSDE \ref{zlogz1} .
We show the value function as a solution of BSDE \ref{zlogz1}.
%%%%%%%%%%%%%%%%%%%%%%%%%%%%%%%%%%%%%%%%%%%%%%%%%%%%%%%%%%%%%%%%%
                 \section{Assumptions and formulation of the problem}

%%%%%%%%%%%%%%%%%%%%%%%%%%%%%%%%%%%%%%%%%%%%%%%%%%%%
%%%%%%%%%%%%%%%%%%%%%%%%%%%%%%%%%%%%%%%%%%%%%%%%%%%%

Let $(\Omega, \mathcal{F}, P)$ be a fixed probability space on which
is defined a standard $d$-dimensional Brownian motion
$B=(B_t)_{0\leq t\leq T}$ whose natural filtration is
$(\mathcal{F}_t^0:=\sigma \{B_s, s\leq t\})_{0\leq t\leq T}$. Let
$\mathbf{F}=(\mathcal{F}_t)_{0\leq t\leq T}$ be
the completed filtration of $(\mathcal{F}_t^0)_{0\leq t\leq T}$ with the $P$%
-null sets of $\mathcal{F}$. We consider the following assumptions,

 %\indent We consider the BSDE,

%\begin{equation}\label{zlogz}
%\displaystyle Y_t = \xi + \int_t^T \varphi(s,Y_s,Z_s) ds - \int_t^T
%Z_s dB_s.
%\end{equation}
%%we put \, \ $f(Z) := |Z|\sqrt{\ln(|Z|)}$

\begin{enumerate}
\item[(\bf {H.1})] \ \ $\E \left[ |\xi|^{\ln(CT+2)+2}\right]<+\infty.$

\item[(\bf{H.2})]
\begin{enumerate}\item[(i)] Assume $\varphi$ is continuous in $(y,z)$ for almost all $(t,w)$;
    \item[(ii)] There exist a constant positive $c_0$ and a process $\eta_t $ satisfying $$\E\left[\int_{0}^{T} \eta_s^{ \ln(Cs+2)+2}ds\right]<+\infty.$$
    and such that for every $t, \omega, y, z$ : \\
    $$\mid \varphi(t,w,y,z)\mid \leq \eta_t+c_0|z|\sqrt{\ln(|z|)}. $$
    \end{enumerate}

\item[(\bf {H.3})]
%There exist $M_1 > 0$ and $1 < \alpha <2$,
% such that:
%$$\mid \varphi(z)\mid \leq  M_1(1 +\mid z\mid^{\alpha}).$$

\par

%(b) \  For every $ z,z'$ such that $ \mid z\mid, \mid z'\mid \leq
%N$, \ \  $\mid \varphi(z)- \varphi(z')\mid \leq
%c\left(\sqrt{\log{N}}\mid z-z'\mid + \dfrac{\log{N}}{N}\right)$ for
%some positive constant $c$ and $N$ large enough.
%\par

  \ There exist $ v\in\L^{q'}(\Omega\times [0, T]; \R_
 +))$ (for some $q'>0$)
 and a real valued sequence $(A_N)_{N>1}$ and constants
 $M_2 \in\R_+$, $r>0$ such that: \\ i) $\forall N>1$, \quad $1<A_N\leq N^{r}.$
\\ ii) $\lim_{N\rightarrow\infty} A_N =
\infty .$
\\ iii) For every $N\in\N,\; \hbox {and every} \ y,\; y'\; z,\; z'
\;\hbox{such that}\; \mid y\mid,\; \mid y'\mid,\; \mid z\mid, \;\mid
z'\mid\leq N$, we have
\begin{align*}
\big(y-y^{\prime}\big)
\big(\varphi(t,\omega,y,z)-\varphi(t,\omega,y^{\prime},z')\big)
\1_{\{v_t(\omega)\leq N\}}
 \leq \
 & M_2\mid
y-y^{\prime}\mid^{2}\log A_{N} \\
& + M_2\mid y-y^{\prime}\mid\mid
z-z^{\prime}\mid\sqrt{\log A_{N}} \\
& + M_2 \dfrac{\log A_{N}}{ A_{N}}.
\end{align*}

%%%%%%%%%%%%%%%%%%%%%%%%%%%%%%%%%%%%%%%%%%%%%%%%%%%%%%%%%%%
%                            proof lemma monotonie
%%%%%%%%%%%%%%%%%%%%%%%%%%%%%%%%%%%%%%%%%%%%%%%%%%%%%%%%%%%

\end{enumerate}

%%%%%%%%%%%%%%%%%%%%%%%%%%%%%%%%%%%%%%%%%%%%%%%%%%%%%%%%%%%%%%%%%%%%%%%%%%%
%                             main result                                 %
%%%%%%%%%%%%%%%%%%%%%%%%%%%%%%%%%%%%%%%%%%%%%%%%%%%%%%%%%%%%%%%%%%%%%%%%%%%
\section{The main results}
 The main objective of this paper is to focus on the
existence and uniqueness of the solution of equation (\ref{zlogz1})
under the previous assumptions.

We denote by $\E$ the set of $\R\times\R^{d}$-valued processes
$(Y,Z)$ defined on $\R_+\times\Omega$ which are $\F_t$-adapted and
such that: $\vert\vert (Y,Z)\vert\vert^2 =
\E\big(\displaystyle\sup_{0\leq t\leq T}\mid Y_t\mid^2 + \int_0^T |
Z_s | ^2ds\big) < +\infty$. The couple $(\E, \vert\vert
.\vert\vert)$ is then a Banach space.

For $N\in \N^*$, we define \begin{equation}\label{rho}
\rho_N(\varphi) = E\displaystyle\int_0^T\sup_{\vert y\vert, \vert
z\vert\leq N}\vert \varphi(s,y,z)\vert ds. \end{equation}
%%%%%%%%%%%%%%%%%%%%%%%%%%%%%%%%%%%%%%%%%%%%%%%%%%%%%%%%%%%%%%%%%%%%%%%%%
%                          definition 1
%%%%%%%%%%%%%%%%%%%%%%%%%%%%%%%%%%%%%%%%%%%%%%%%%%%%%%%%%%%%%%%%%%%%%%%%
\begin{definition}
 A solution of equation
(\ref{zlogz1}) is a couple $(Y,Z)$ which belongs to the space $(\E,
\vert\vert.\vert\vert)$ and satisfies equation (\ref{zlogz1}).
\end{definition}
%%%%%%%%%%%%%%%%%%%%%%%%%%%%%%%%%%%%%%%%%%%%%%%%%%%%%%%%%%%%%%%%%%%%%%%%%
%                          Assumptions
%%%%%%%%%%%%%%%%%%%%%%%%%%%%%%%%%%%%%%%%%%%%%%%%%%%%%%%%%%%%%%%%%%%%%%%%

%%%%%%%%%%%%%%%%%%%%%%%%%%%%%%%%%%%%%%%%%%%%%%%%%%%%%%%%%%%%%%%%%%%%%%%%%%%
%                               Theorem uniqueness
%%%%%%%%%%%%%%%%%%%%%%%%%%%%%%%%%%%%%%%%%%%%%%%%%%%%%%%%%%%%%%%%%%%%%%%%%%%

The main result of this section are the following two theorems.
\begin{theorem}\label{unique}
  Assume that
\bf{(H.1)},\bf{(H.2)} and \bf{(H.3)} are satisfied. Then, equation
(\ref{zlogz1}) has a unique solution.
\end{theorem}

%%%%%%%%%%%%%%%%%%%%%%%%%%%%%%%%%%%%%%%%%%%%%%%%%%%%%%%%%%%
%                            Stability
%%%%%%%%%%%%%%%%%%%%%%%%%%%%%%%%%%%%%%%%%%%%%%%%%%%%%%%%%%%
In the following, we give a stability result for the solution with
respect to the data $(\varphi,\xi)$. Roughly speaking, if
$\varphi_n$ converges to $\varphi$ in the metric defined by the
family of semi-norms $(\rho_N)$ and $\xi_n$ converges to $\xi$ in
$L^{2}(\Omega)$ then $(Y^n,Z^n)$ converges to $(Y,Z)$ in some
reflexive Banach space which we will precise below. Let
$(\varphi_n)$ be a sequence of functions which are measurable for
each $n$. Let $(\xi_n)$ be a sequence of random variables which are
$\F_T$-measurable for each $n$ and such that \ $\displaystyle
sup_n\E \big( |\xi_n|^{\ln(CT+2)+2}\big)<+\infty$. We will assume
that for each $n$, the BSDE
%$(E^{\varphi_n,\xi_n})$
corresponding to the data $(\varphi_n,\xi_n)$ has a (not necessarily
unique) solution. Each solution of the BSDE $(\varphi_n,\xi_n)$ will
be denoted by $(Y^n,Z^n)$. We consider the following assumptions,

\vskip 0.2cm\noindent \bf{(H.5)} \quad {\it For every $N$, \  $
\rho_N(\varphi_n-\varphi) \longrightarrow 0 $ as $n \longrightarrow
\infty$.
\par\noindent
\bf{(H.6)} \quad $E(\vert \xi_n - \xi \vert^{\ln(CT+2)+2})\longrightarrow 0$
as $n\rightarrow \infty$ .
\\
%\bf{(H.7)}\; There exist $M_1 > 0$, $1< \alpha <2$, $\alpha'>1$ and
%$\overline{\eta}\in {\L}^{\alpha'}([0,T]\times\Omega)$ such that:
%$$\displaystyle\sup_{n}\mid \varphi_n(t,\omega, y, z)\mid \leq \overline{\eta}
%+M_1(\mid y\mid^{\alpha}+\mid z\mid^{\alpha}).$$
%
%\par\noindent
\bf{(H.7)}\; there exist a constant positive $c_0$ and $\eta_t $
satisfying
$$\E\left[\int_{0}^{T} \eta_s^{ \ln(Cs+2)+2}ds\right]<+\infty,$$
    and such that:\\
    $$\displaystyle\sup_{n}\mid \varphi_n(t,\omega, y, z)\mid \leq \eta_t+c_0|z|\sqrt{\ln(|z|)}$$

%%%%%%%%%%%%%%%%%%%%%%%%%%%%%%%%%%%%%%%%%%%%%%%%%%%%%%%%%%%
%                            Theorem 2.2
%%%%%%%%%%%%%%%%%%%%%%%%%%%%%%%%%%%%%%%%%%%%%%%%%%%%%%%%%%%
\begin{theorem}\label{stability}  Let $\varphi$ and $\xi$
be as in Theorem \ref{unique}. Assume that \bf{(H.5)}, \bf{(H.6)},
and \bf{(H.7)} are satisfied. Then, for all $q<2$  we
have $$ \lim_{n\rightarrow +\infty} \left(\E\sup_{0\leq t\leq
T}\vert Y_{t}^n - Y_{t}\vert^{q}+\E\dint_{0}^{T}\vert Z_{s}^n -
Z_{s}\vert^{q}ds\right) = 0.
$$
\end{theorem}

%%%%%%%%%%%%%%%%%%%%%%%%%%%%%%%%%%%%%%%%%%%%%%%%%%%%%%%%%%%
%                            remark y alpha
%%%%%%%%%%%%%%%%%%%%%%%%%%%%%%%%%%%%%%%%%%%%%%%%%%%%%%%%%%%
\begin{remark}\label{remark}   The conclusions of the previous theorems remain valid if, instead of hypothesis \bf{(H2)-(ii)}, we assume the following more general condition :
\end{remark}
\noindent\bf{(H2)-(iii)} \ There exist a constants positive $c_0$, $0 < \alpha'< 2$  and a process $\eta_t $ satisfying $$\E\left[\int_{0}^{T} \eta_s^{ \ln(Cs+2)+2}ds\right]<+\infty,$$
    and such that for every $t, \omega, y, z$:\\
    $$\mid \varphi(t,w,y,z)\mid \leq \eta_t+ \vert y\vert^{\alpha'} + c_0|z|\sqrt{\ln(|z|)}$$

%%%%%%%%%%%%%%%%%%%%%%%%%%%%%%%%%%%%%%%%%%%%%%%%%%%%%%%%%%%%%%%%%%%%%%%%%%%%%%%%
%                          Proofs
%%%%%%%%%%%%%%%%%%%%%%%%%%%%%%%%%%%%%%%%%%%%%%%%%%%%%%%%%%%%%%%%%%%%%%%%%%%%%%%%
\section{Proofs}

To prove Theorem \ref{unique} and Theorem \ref{stability}, we need
the following lemmas.

\begin{lemma} \label{estimateY}
Let $(Y,Z)$ be a solution of the above BSDE, where $(\xi,\varphi)$
satisfies the assumptions $(\bf {H1})$ and $(\bf {H2})$. Then there
exists a constant $C_{T},$ such that:
$$\E \sup_{t \in [0,T]} |Y_t|^{\ln(Ct+2)+2} \leq C_{T}\E \left[ |\xi|^{\ln(CT+2)+2} + \int_0^{T}\eta_s^{ \ln(Cs+2)+2}ds \right ].$$
\end{lemma}

%%%%%%%%%%%%%%%%%%%%%%%%%%%%%%%%%%%%%%%%%%%%%%%%%%%%%%%%%%%
%                            Proof estimateY
%%%%%%%%%%%%%%%%%%%%%%%%%%%%%%%%%%%%%%%%%%%%%%%%%%%%%%%%%%%

\bop. For some constant $C$ large, let us consider the function from
$[0,T]\times\R$ into $\R^+$ defined by.
$$u(t,x)=\mid x\mid^{\ln(Ct+2)+2}.$$
Then

$$u_t=\frac{C}{Ct+2}\ln(\mid x\mid)\mid x\mid^{\ln(Ct+2)+2},\qquad u_x=(\ln(Ct+2)+2)\mid x\mid^{\ln(Ct+2)+1} sgn(x)$$
and $u_{xx}=(\ln(Ct+2)+2)(\ln(Ct+2)+1)\mid x\mid^{\ln(Ct+2)}$, with
the notation $sgn(x)=-\ind_{x\leq0}+\ind_{x>0}$. For $k\geq0$, let
$\tau_{k}$ be the stopping time defined as follows:
$$\tau_k=\inf \{t\geq 0,\, \int_{0}^{T}(\ln(Cs+2)+2)^2\mid
Y_s\mid^{2\ln(Cs+2)+2}\mid Z_s\mid^2 ds]\geq k\}\wedge T.$$

Next using It\^o's formula yields:
\begin{eqnarray*}
&& \mid Y_{t\wedge \tau_k}\mid^{\ln(Ct+2)+2} =\mid
Y_{\tau_k}\mid^{\ln(Ct+2)+2}- \int_{t\wedge \tau_k}^{\tau_k}\frac{C}{Cs+2}\ln(\mid Y_s \mid)\mid  Y_s\mid^{\ln(Cs+2)+2}ds \\
\nonumber && \qquad \qquad -\frac{1}{2}\int_{t\wedge
\tau_k}^{\tau_k} |Z_s|^2(\ln(Cs+2)+2)(\ln(Cs+2)+1)\mid
Y_s\mid^{\ln(Cs+2)}ds
\\\nonumber && \qquad \qquad +\int_{t\wedge \tau_k}^{\tau_k}
(\ln(Cs+2)+2)\mid Y_s\mid^{\ln(Cs+2)+1}sgn(Y_s)f(s,Y_s,Z_s)ds  \\
\nonumber && \qquad \qquad -  \int_{t\wedge
\tau_k}^{\tau_k}(\ln(Cs+2)+2)\mid
Y_s\mid^{\ln(Cs+2)+1}sgn(Y_s) Z_s dB_s, \\
\nonumber && \qquad \qquad \leq \mid
Y_{\tau_k}\mid^{\ln(Ct+2)+2}- \int_{t\wedge \tau_k}^{\tau_k}\frac{C}{Cs+2}\ln(\mid Y_s \mid)\mid  Y_s\mid^{\ln(Cs+2)+2}ds \\
\nonumber && \qquad \qquad -\frac{1}{2}\int_{t\wedge
\tau_k}^{\tau_k} |Z_s|^2(\ln(Cs+2)+2)(\ln(Cs+2)+1)\mid
Y_s\mid^{\ln(Cs+2)}ds
\\\nonumber && \qquad \qquad +\int_{t\wedge \tau_k}^{\tau_k}
(\ln(Cs+2)+2)\mid Y_s\mid^{\ln(Cs+2)+1}(\eta_s+c_0|Z_s|\sqrt{\ln(|Z_s|)})ds  \\
\nonumber && \qquad \qquad -  \int_{t\wedge
\tau_k}^{\tau_k}(\ln(Cs+2)+2)\mid Y_s\mid^{\ln(Cs+2)+1}sgn(Y_s) Z_s
dB_s.
\end{eqnarray*}
By Young's inequality it hold true that:
$$
(\ln(Cs+2)+2)\mid Y_s\mid^{\ln(Cs+2)+1}\eta_s \leq \mid
Y_s\mid^{\ln(Cs+2)+2}+ (\ln(Cs+2)+2)^{\ln(Cs+2)+1}\eta_s^{
\ln(Cs+2)+2}.
$$
For $\mid y \mid$ large enough and the last inequality there exists
$C_1$ such that:
\begin{eqnarray*}
&& \mid Y_{t\wedge \tau_k}\mid^{\ln(Ct+2)+2} =\mid
Y_{\tau_k}\mid^{\ln(Ct+2)+2}- \int_{t\wedge \tau_k}^{\tau_k}C_1\ln(\mid Y_s \mid)\mid  Y_s\mid^{\ln(Cs+2)+2}ds \\
\nonumber && \qquad \qquad -\frac{1}{2}\int_{t\wedge
\tau_k}^{\tau_k} |Z_s|^2(\ln(Cs+2)+2)(\ln(Cs+2)+1)\mid
Y_s\mid^{\ln(Cs+2)}ds
\\\nonumber && \qquad \qquad +\int_{t\wedge \tau_k}^{\tau_k}
(\ln(Cs+2)+2)\mid Y_s\mid^{\ln(Cs+2)+1}c_0|Z_s|\sqrt{\ln(|Z_s|)}ds
\\\nonumber && \qquad \qquad
 +\int_{t\wedge \tau_k}^{\tau_k}(\ln(Cs+2)+2)^{\ln(Cs+2)+1}\eta_s^{ \ln(Cs+2)+2}ds \\
\nonumber && \qquad \qquad -  \int_{t\wedge
\tau_k}^{\tau_k}(\ln(Cs+2)+2)\mid
Y_s\mid^{\ln(Cs+2)+1}sgn(Y_s) Z_s dB_s, \\
\nonumber && \qquad \qquad \leq \mid
Y_{\tau_k}\mid^{\ln(Ct+2)+2}- \int_{t\wedge \tau_k}^{\tau_k}(\ln(Cs+2)+2)(\ln(Cs+2)+1)\mid Y_s\mid^{\ln(Cs+2)}[ \\
\nonumber && \qquad \qquad \frac{C_1\ln(\mid Y_s \mid)\mid
Y_s\mid^{2}}{(\ln(Cs+2)+2)(\ln(Cs+2)+1)}+\frac{|Z_s|^2}{2}-\frac{(\ln(Cs+2)+2)\mid
Y_s\mid c_0|Z_s|\sqrt{\ln(|Z_s|)}}{(\ln(Cs+2)+2)(\ln(Cs+2)+1)}]ds
\\\nonumber && \qquad \qquad +\int_{t\wedge \tau_k}^{\tau_k}(\ln(Cs+2)+2)^{\ln(Cs+2)+1}\eta_s^{ \ln(Cs+2)+2}ds\\
\nonumber && \qquad \qquad -  \int_{t\wedge
\tau_k}^{\tau_k}(\ln(Cs+2)+2)\mid Y_s\mid^{\ln(Cs+2)+1}sgn(Y_s) Z_s
dB_s.
\end{eqnarray*}
There exist a constants $C_2$ and $C_3$ ($C_2>2C_3^2$)

\begin{eqnarray}\label{fin-est}
&& \mid Y_{t\wedge \tau_k}\mid^{\ln(Ct+2)+2}  \leq
\\\nonumber && \qquad \qquad \mid
Y_{\tau_k}\mid^{\ln(Ct+2)+2}- \int_{t\wedge \tau_k}^{\tau_k}(\ln(Cs+2)+2)(\ln(Cs+2)+1)\\
\nonumber && \qquad \qquad \mid Y_s\mid^{\ln(Cs+2)}[ C_2\ln(\mid Y_s
\mid)\mid Y_s\mid^{2}+\frac{|Z_s|^2}{2}-C_3\mid
Y_s\mid|Z_s|\sqrt{\ln(|Z_s|)}]ds
\\\nonumber && \qquad \qquad +\int_{t\wedge \tau_k}^{\tau_k}(\ln(Cs+2)+2)^{\ln(Cs+2)+1}\eta_s^{ \ln(Cs+2)+2}ds\\
\nonumber && \qquad \qquad -  \int_{t\wedge
\tau_k}^{\tau_k}(\ln(Cs+2)+2)\mid Y_s\mid^{\ln(Cs+2)+1}sgn(Y_s) Z_s
dB_s.
\end{eqnarray}

Now we show that \begin{equation}\label{majoration}C_3\mid
Y_s\mid|Z_s|\sqrt{\ln(|Z_s|)}\leq \frac{|Z_s|^2}{2}+ C_2\ln(\mid
Y_s\mid)\mid Y_s\mid^{2}\end{equation}

if $|Z_s|\leq \mid Y_s\mid$, (\ref{majoration}) is obviously true.
Assume $|Z_s|> \mid Y_s\mid$. Denote $a_s=\frac{|Z_s|}{|Y_s|}$.
Then\\
$C_3\mid Y_s\mid|Z_s|\sqrt{\ln(|Z_s|)}\leq C_3 a_s
Y_s^2[\sqrt{\ln(|Z_s|)}+\sqrt{\ln(|Y_s|)}]$,\\
$\frac{|Z_s|^2}{2}+ C_2\ln(\mid Y_s\mid)\mid Y_s\mid^{2}\leq
[\frac{a_s^2}{2}+ C_2\ln(\mid Y_s\mid)]\mid Y_s\mid^{2}$.\\
Obviously $$C_3 a_s [\sqrt{\ln(|Y_s|)}\leq
\frac{1}{2}[\frac{a_s^2}{2}+2C_3^2\ln(\mid Y_s\mid)].$$ Assume $r$
is the constant such that $C_3\sqrt{\ln(r)}=\frac{r}{4}$. If
$a_s\geq r,$
$$C_3a_s\sqrt{\ln(a_s)}\leq \frac{a_s^2}{4}.$$
If $a_s\leq r,$ and $\mid y \mid$ large enough then
$$C_3a_s\sqrt{\ln(a_s)}\leq C_3r\sqrt{\ln(r)}\leq \frac{C_2}{2}\ln(\mid Y_s\mid).$$
Then (\ref{majoration}) holds. Finally taking the limit in both
sides as $k\rightarrow + \infty$ and the lemma is proved. \eop

%%%%%%%%%%%%%%%%%%%%%%%%%%%%%%%%%%%%%%%%%%%%%%%%%%%%%%%%%%%
%                            lemma estimate Z
%%%%%%%%%%%%%%%%%%%%%%%%%%%%%%%%%%%%%%%%%%%%%%%%%%%%%%%%%%%

\bigskip
\begin{lemma} \label{estimateZ}
Let $(Y,Z)$ be a solution of the above BSDE. Then
 There exits a real constant $C_{p}$ depending only on $p$
 such that:
$$\E \left[ \left( \int_0^T |Z_s|^2 ds \right)^{p/2} \right] \leq C_{p} \E \left[ |\xi|^p+\sup_{t \in [0,T]} |Y_t|^{p\frac{2+\ln(2)}{2}} + \left(\int_0^T\mid
\eta_s \mid^2ds\right)^{\frac{p}{2}} \right].$$
\end{lemma}

\bop.  Applying It\^o's formula to the process $Y_t$ and the
function $y\longmapsto y^2$ yields:
\begin{eqnarray*}
&& |Y_0|^2 + \int_0^{T} |Z_s|^2 ds = |\xi|^2 + 2\int_0^{T} Y_s
\varphi(s,Y_s,Z_s) ds - 2 \int_0^{T}  Y_s Z_s dB_s \\
&& \leq |\xi|^2 + 2\int_0^{T}\mid Y_s\mid(\mid
\eta_s\mid+c_0|Z_s|\sqrt{\ln(|Z_s|)}) ds\mid - 2 \int_0^{T} Y_s Z_s
dB_s .\end{eqnarray*} As we have $$2\mid Y_s\mid \mid \eta_s\mid\leq
\mid Y_s\mid^2+\mid \eta_s\mid^2,$$ and for any $\varepsilon>0$ we
have:

 $$\sqrt{2\varepsilon \ln(|z|)}=\sqrt{\ln(|z|^{2\varepsilon })}\leq
|z|^{\varepsilon }.$$ Then plug the two last inequalities in the
previous one to obtain:
\begin{eqnarray*}\\
&& |Y_0|^2 + \int_0^{T} |Z_s|^2 ds \\&&\leq |\xi|^2 + \sup_{s\leq
T}|Y_s^2|+\int_0^T\mid \eta_s
\mid^2ds+\frac{2}{\sqrt{2\varepsilon}}\int_0^T \mid Y_s\mid \mid
Z_s\mid^{1+\varepsilon} ds|)- 2 \int_0^{T} Y_s Z_s dB_s.
\end{eqnarray*}
We now choose $0<\varepsilon <1$ and by young's inequality it holds
true that:
$$2\frac{\mid Y_s\mid}{\sqrt{2\varepsilon}} \mid
Z_s\mid^{1+\varepsilon} \leq \frac{ 1-\varepsilon
}{2}(\frac{2}{\sqrt{2\varepsilon}})^{\frac{2}{1-\varepsilon}}\mid Y_s\mid
^{\frac{2}{1-\varepsilon}}+\frac{1+\varepsilon}{2}\mid Z_s\mid^{2}.$$
Then, there exists a positive constant $c_\varepsilon$

\begin{eqnarray*}
&& |Y_0|^2 + \int_0^{T} |Z_s|^2 ds \leq |\xi|^2 + \sup_{s\leq
T}|Y_s^2|+\int_0^T\mid \eta_s \mid^2ds\\&&
\qquad\qquad+c_\varepsilon\sup\limits_{s\leq T}\mid Y_s\mid
^{\frac{2}{1-\varepsilon}}+\frac{1+\varepsilon}{2} \int_0^T\mid
Z_s\mid^{2}ds- 2 \int_0^{T} Y_s Z_s dB_s.
\end{eqnarray*}

For $\mid y \mid$ large enough and $\varepsilon \leq
\frac{\ln(2)}{2+\ln(2)}$ then

\begin{eqnarray*}
&& |Y_0|^2 + \int_0^{T} |Z_s|^2 ds \leq |\xi|^2 + c_\varepsilon
\sup_{s\leq T}|Y_s|^{2+\ln(2)}+\int_0^T\mid \eta_s \mid^2ds\\&&
\qquad\qquad\qquad+\frac{1+\varepsilon}{2} \int_0^T\mid Z_s\mid^{2}ds-
2 \int_0^{T} Y_s Z_s dB_s.
\end{eqnarray*}

Then we obtain:
\begin{eqnarray*}
&&\E \left( \int_0^{T} |Z_s|^2 ds \right)^{p/2} \leq
C_p\E\left[c_\varepsilon\sup_{s\leq T}|Y_s|^{2+\ln(2)}+\int_0^T\mid
\eta_s \mid^2ds+\left(\frac{1+\varepsilon}{2}\right)^{\frac{p}{2}}
\left(\int_0^T\mid Z_s\mid^{2}ds\right)^{\frac{p}{2}}\right]\\&&
\qquad\qquad\qquad\qquad\qquad\qquad+ C_p\E\left[\left|\int_0^{T}
Y_s Z_s dB_s\right|^{\frac{p}{2}}\right].
\end{eqnarray*}
Next thanks to BDG's inequality and for any $\beta >0$ we have:
\begin{eqnarray*}
&&\E \left[  \left| \int_0^{t}Y_s Z_s dB_s\right|^{p/2} \right] \leq
 \bar C_p \E \left[  \left( \int_0^{T} |Y_s|^2 |Z_s|^2 ds
\right)^{p/4}\right]  \\
&& \leq  \bar C_p (\E \left[  \left( \sup_{t \in [0,T]} |Y_t|
\right)^{p/2}  \left( \int_0^{T} |Z_s|^2 ds \right)^{p/4}
+\eps^{\frac{p}{2}}\left(\int_0^T |Z_s|^2
ds\right)^{p/2}\right] \\
&& \leq  \frac{\bar C_p^2}{\beta} \E \left[ \sup_{t \in [0,T]}
|Y_t|^p \right] + \beta \E\left[  \left( \int_0^{T} |Z_s|^2 ds
\right)^{p/2}\right].
\end{eqnarray*} Choosing
$\beta$ and $\eps$ small enough to obtain the desired result. \eop

\begin{lemma}\label{estimatevarphi} If \bf{(H.2)}
 holds  then,
\begin{align*}
&  \E  \dint_0^T   | \varphi(s,Y_s,Z_s)|^{\overline{\alpha}} ds  \; \leq \;
K\big[1 + \E\dint_0^T
 {\eta}_s^2 ds+
\E\dint_0^T\vert Z_s\vert^2 ds\big]
\end{align*}
where $\overline{\alpha}=\min(2,
\dfrac{2}{\alpha})$ and $K$ is a positive constant which depends on $c_0$ and $T$.
\end{lemma}

%%%%%%%%%%%%%%%%%%%%%%%%%%%%%%%%%%%%%%%%%%%%%%%%%%%%%%%%%%%%%%%%%%%%
%                       Proof of Prop 3.2                  %
%%%%%%%%%%%%%%%%%%%%%%%%%%%%%%%%%%%%%%%%%%%%%%%%%%%%%%%%%%%%%%%%%%%

\noindent\bf{Proof.} Observe that assumption (\bf{H.2}) implies that there exist $c_1 > 0$ and $0\leq \alpha <2$
 such that:
 \begin{equation}\label{aplpha}
\mid\varphi(t,\omega, y, z)\mid \leq {\eta}_t + c_1\mid z\mid^{\alpha}.
\end{equation}
We successively use Assumption (H.3) and inequality (\ref{aplpha}) to show that
\begin{align*}
 \E  \dint_0^T   | f(s,Y_s,Z_s)|^{\overline{\alpha}} ds
 &\leq
\E\dint_0^T (\eta_s+c_0|z|\sqrt{\ln(|z|)})^{\overline{\alpha}} ds
\\ &\leq \E\dint_0^T
({\eta}_s^{} + c_1\vert
Z_s\vert^{\alpha})^{\overline{\alpha}} ds
\\ &\leq
(1 + c_1^{\overline{\alpha}})\E\dint_0^T (( {\eta}_s)^{\overline{\alpha}} + ( \vert Z_s\vert)^{\alpha\overline{\alpha}}) ds\\
&\leq (1 + c_1^{\overline{\alpha}})\E\dint_0^T (( 1+ {\eta}_s)^{\overline{\alpha}}  + (1 + \vert Z_s\vert)^{\alpha\overline{\alpha}}) ds \\
&\leq (1 + c_1^{\overline{\alpha}})\E\dint_0^T (( 1+ {\eta}_s)^{2}  + (1 + \vert Z_s\vert)^{2}) ds
\\ &\leq (1 + c_1^{\overline{\alpha}})\big(4T + \E\dint_0^T
( {\eta}_s^{2} + \vert Z_s\vert^{2})
ds\big)
\end{align*}
Lemma \ref{estimatevarphi} is proved.
\eop

%%%%%%%%%%%%%%%%%%%%%%%%%%%%%%%%%%%%%%%%%%%%%%%%%%%%%%%%%%%%%%%%%%%%
%                         Lemma 3.1
%%%%%%%%%%%%%%%%%%%%%%%%%%%%%%%%%%%%%%%%%%%%%%%%%%%%%%%%%%%%%%%%%%%%
\begin{lemma}\label{lem1}
{\it There exists a sequence of functions $(\varphi_n)$ such that,
\par\noindent
$(a)$ \; For each $n$, $\varphi_n$ is bounded and globally Lipschitz
in $(y,z)$ $a.e.$ $t$ and $P$-$a.s.\omega$.
\\
$(b)$ \; $\displaystyle\sup_{n}\mid \varphi_n(t,\omega, y, z)\mid
\leq \eta_t+c_0|z|\sqrt{\ln(|z|)}$, \quad $P$-$a.s.$, $a.e.$ $t\in
[0,T]$.
%\\
%$\sup_n\vert \varphi_n(t,\omega,y,z)\vert\; \leq M' +M_1(\mid
%y\mid^{\alpha}+\mid z\mid^{\alpha})$.\quad $P$-$a.s.$, $a.e.$ $t\in
%[0,T]$.
\\
 $(c)$\; For every $N$,
$\rho_N (\varphi_n-\varphi)\longrightarrow 0$ as
$n\longrightarrow\infty$. }
\end{lemma}
%%%%%%%%%%%%%%%%%%%%%%%%%%%%%%%%%%%%%%%%%%%%%%%%%%%%%%%%%%%%%%%%%%%%
%                        Proof of Lemma 3.1
%%%%%%%%%%%%%%%%%%%%%%%%%%%%%%%%%%%%%%%%%%%%%%%%%%%%%%%%%%%%%%%%%%%%
\bop Let $\varepsilon_n: \R^2 \longrightarrow \R_+$ be a sequence of
smooth functions with compact support which approximate the Dirac
measure at 0 and which satisfy $\int \varepsilon_n (u)du = 1$. Let
$\psi_n$ from $\R^{2}$ to $ \R_+$  be a sequence of smooth functions
such that $0\leq \vert\psi_n\vert \leq 1$, $\psi_n(u)=1$ for $\vert
u\vert \leq n$ and $\psi_n(u)=0$ for $\vert u\vert \geq n+1$.  We
put, $\varepsilon_{q,n}(t,y,z) = \int
\varphi(t,(y,z)-u)\alpha_q(u)du\psi_n(y,z)$. For $n \in \N^*$, let
$q(n)$ be an integer such that $q(n) \geq n+n^\alpha$. It is not
difficult to see that the sequence $\varphi_n :=
\varepsilon_{q(n),n}$ satisfies all the assertions $(a)$-$(c)$. \eop

%%%%%%%%%%%%%%%%%%%%%%%%%%%%%%%%%%%%%%%%%%%%%%%%%%%%%%%%%%%%%%%%%%%%
%                            Lemma 3.2                             %
%%%%%%%%%%%%%%%%%%%%%%%%%%%%%%%%%%%%%%%%%%%%%%%%%%%%%%%%%%%%%%%%%%%%
Using Lemma \ref{estimateY}, Lemma \ref{estimateZ}, Lemma
\ref{estimatevarphi}, Lemma \ref{lem1} and standard arguments of
BSDEs, one can prove the following estimates.
\begin{lemma}\label{lem2}
Let $\varphi$ and $\xi$ be as in Theorem \ref{unique}. Let
$(\varphi_n)$ be the sequence of functions associated to $\varphi$
by Lemma \ref{lem1}. Denote by $(Y^{\varphi_n},Z^{\varphi_n})$ the
solution of equation $(E^{\varphi_n})$. Then, there exit
 constants $K_1$, $K_2$, $K_3$  and a universal constant $\ell$ such that
\\$a)$\, $ \displaystyle\sup_n\E\int_0^T \vert Z_s^{\varphi_n}\vert^2ds
\leq
K_1 $ \\
$b)$\, $ \displaystyle \sup_n\E\sup_{0\leq t\leq T}(\mid
Y_t^{\varphi_n}\mid^2) \leq \ell K_1 :=
 K_2$
\\$c)$\, $\displaystyle \sup_n\E\dint_0^T \vert \varphi_n(s,Y_{s}^{\varphi_n},
Z_{s}^{\varphi_n})\vert^{\overline{\alpha}}ds \leq
 K_{3}$
 \\
 where $\overline{\alpha}=\min(2,
\dfrac{2}{\alpha})$
\end{lemma}

%%%%%%%%%%%%%%%%%%%%%%%%%%%%%%%%%%%%%%%%%%%%%%%%%%%%%%%%%%%%%%%%%%%%
%                        Proof of Lemma 3.2
%%%%%%%%%%%%%%%%%%%%%%%%%%%%%%%%%%%%%%%%%%%%%%%%%%%%%%%%%%%%%%%%%%%%
%\bop \ Using It\^o's formula and lemma 3.1 (c), we show that for all
%$t\leq T$ \vskip 0.2cm $ \displaystyle e^{2Mt}\mid
%Y_t^{\varphi_n}\mid^2 + (1-2\gamma)\dint_t^Te^{2Ms}\mid
%Z_s^{\varphi_n}\mid^2ds \leq
% e^{2MT}\mid \xi\mid^2 +
%2\dint_t^Te^{2Ms}(\eta_s+M')ds
%$
%\par\hskip 7.4cm \quad
%$ \displaystyle -2\dint_t^Te^{2Ms}\langle Y_s^{\varphi_n},\quad
%Z_s^{\varphi_n}dW_s\rangle $
%\par
%Taking expectation we get assertion $a)$. Assertion $b)$ is a
%direct consequence of the Burkholder Davis Gundy inequality and assertion
%$a)$. Finally, assertions $c)$ and $d)$  follow from Lemma 3.4
%$(b)$ and assumption (H.3). Lemma 3.5 is proved.
% \eop
%%%%%%%%%%%%%%%%%%%%%%%%%%%%%%%%%%%%%%%%%%%%%%%%%%%%%%%%%%%%%%%%%%%%
%                       corollary 3.1                          %
%%%%%%%%%%%%%%%%%%%%%%%%%%%%%%%%%%%%%%%%%%%%%%%%%%%%%%%%%%%%%%%%%%%
After extracting a subsequence, if necessary, we have
\begin{corollary}
There are $Y\in\L^{2}(\Omega, L^{\infty}[0,T])$,
$Z\in\L^{2}(\Omega\times[0,T])$,
$\Gamma\in\L^{\overline{\alpha}}(\Omega\times[0,T])$ such that
\begin{align*}
& Y^{\varphi_n}\rightharpoonup Y,\; \hbox{weakly star in}\;
\L^{2}(\Omega, L^{\infty}[0,T])
\\ & Z^{\varphi_n}\rightharpoonup Z,\;
\hbox{weakly in}\; \L^{2}(\Omega\times[0,T])
\\ & \varphi_n(., Y^{\varphi_n}, Z^{\varphi_n}) \rightharpoonup \Gamma_{.} \;
\hbox{weakly in}\; \L^{\overline{\alpha}}(\Omega\times[0,T]),
\end{align*}
and moreover
\begin{align*}
Y_{t}= \xi + \dint_{t}^{T}\Gamma_{s}ds-\dint_{t}^{T}Z_{s}dW_{s},
\; \forall t\in[0,T].
\end{align*}
\end{corollary}

%%%%%%%%%%%%%%%%%%%%%%%%%%%%%%%%%%%%%%%%%%%%%%%%%%%%%%%%%%%%%%%%%%%%%%%%
%                       lemma 3.3                                      %
%%%%%%%%%%%%%%%%%%%%%%%%%%%%%%%%%%%%%%%%%%%%%%%%%%%%%%%%%%%%%%%%%%%%%%%%

The following lemma, were established in \cite{BEHP}, is a direct
consequence of H\"{o}lder's and Schwarz's inequalities and the fact that $ab\leq \dfrac{\alpha^{2}}{2}a^2
+\dfrac{1}{2\alpha^{2}}b^2$ \ for each $\alpha >0 $ and each real numbers $a$, $b$.
\begin{lemma}\label{lem3}
For every $\beta \in ]1,2]$, $A>0$, $(y)_{i=1..d}\subset \R$,
$(z)_{i=1..d,j=1..r}\subset \R$ we have,
\begin{align*}
 A\vert y\vert\vert z\vert - \dfrac{1}{2}\vert z\vert^2
 +\dfrac{2-\beta}{2}\vert y\vert^{-2}\vert y z\vert^2
 \leq
 \dfrac{1}{\beta -1}A^2\vert y\vert^2 - \dfrac{\beta -1}{4}\vert z\vert^2.
\end{align*}
\end{lemma}

This lemma remains valid in multidimensional case.
%%%%%%%%%%%%%%%%%%%%%%%%%%%%%%%%%%%%%%%%%%%%%%%%%%%%%%%%%%%%%%%%%%%%%%%%
%                       Proof of lemma 3.3                                      %
%%%%%%%%%%%%%%%%%%%%%%%%%%%%%%%%%%%%%%%%%%%%%%%%%%%%%%%%%%%%%%%%%%%%%%%%

%\bop  Using the inequality $ab\leq \dfrac{\alpha^{2}}{2}a^2
%+\dfrac{1}{2\alpha^{2}}b^2$ we have
%\begin{align*}
% A\vert y\vert\vert z\vert - \dfrac{1}{2}\vert z\vert^2
% +\dfrac{2-\beta}{2}\dfrac{1}{\vert y\vert^{2}}\vert y z\vert^2
% \leq
% \dfrac{\alpha^2}{2}A^2\vert y\vert^2 + \dfrac{1}{2\alpha^2}\vert z\vert^2
% -  \dfrac{1}{2}\vert z\vert^2 +
% \dfrac{2-\beta}{2}\dfrac{1}{\vert y\vert^{2}}\vert y z\vert^2.
%\end{align*}
%The proof is finished by choosing $\alpha^{2}=\dfrac{2}{(\beta-1)}$.
%Lemma \ref{lem3} is proved. \eop

%%%%%%%%%%%%%%%%%%%%%%%%%%%%%%%%%%%%%%%%%%%%%%%%%%%%%%%%%%%%%%%%%%%%
%                        estimate between two solutions                     %
%%%%%%%%%%%%%%%%%%%%%%%%%%%%%%%%%%%%%%%%%%%%%%%%%%%%%%%%%%%%%%%%%%%
\vskip 0.3cm
\section {Estimate between two solutions}
The key estimate is given by,
\begin{lemma}\label{lem4}
For every $R\in\N$, $\beta \in
]1,\min\left(3-\frac{2}{\overline\alpha},2\right)[$,
 $\delta' <
 (\beta-1)\min\left(\frac{1}{4M_2^2},\frac{3-\frac{2}{\overline\alpha}-\beta}{2rM_2^2\beta}\right)$ and
 $\varepsilon>0$, there exists $N_{0}>R$ such that for all
 $N>N_{0}$ and $T'\leq T$:
\begin{align*}
& \limsup_{n,m\rightarrow +\infty}E \sup_{(T'-\delta')^{+}\leq t
\leq T'}\vert Y_{t}^{\varphi_n}-Y_{t}^{\varphi_m}\vert^\beta  + E
\int_{(T'-\delta')^{+}}^{T'}{\left|
Z_{s}^{\varphi_n}-Z_{s}^{\varphi_m}\right|^{2}\over\left(\vert
Y_{s}^{\varphi_n}-Y_{s}^{\varphi_m}\vert^{2}+\nu_{R}\right)^{{2-\beta\over
2}}}ds
\\
&\qquad\qquad\qquad\leq \varepsilon +\frac{\ell}{\beta -1}
e^{C_N\delta'} \limsup_{n,m\rightarrow +\infty}E \vert
Y_{T'}^{\varphi_n}-Y_{T'}^{\varphi_m}\vert^\beta.
\end{align*}
where $\nu_{R} = \sup \left\{(A_N)^{-1}, N\geq R
\right\}$,
 $C_N = {2M_2^2\beta\over (\beta -1)} \log A_{N}$ and $\ell$ is a
 universal positive constant.
\end{lemma}

%%%%%%%%%%%%%%%%%%%%%%%%%%%%%%%%%%%%%%%%%%%%%%%%%%%%%%%%%%%%%%%%%%%%
%                       Proof of Lemma 3.4                      %
%%%%%%%%%%%%%%%%%%%%%%%%%%%%%%%%%%%%%%%%%%%%%%%%%%%%%%%%%%%%%%%%%%%

\bop. To simplify the computations, we assume (without loss of
generality) that assumption \bf{(H3)}-C)-iii) holds without the
multiplicative term $\1_{\{v_t(\omega)\leq N\}}$.
\par\noindent
Let $0<T'\leq T$. It follows from It\^o's formula that for all
$t\leq T'$,
\begin{align*}
&\left| Y_{t}^{\varphi_n}-Y_{t}^{\varphi_m}\right|
^{2}+\int_{t}^{T'}\left| Z_{s}^{\varphi_n}-Z_{s}^{\varphi_m}\right|
^{2}ds
\\& = \left| Y_{T'}^{\varphi_n}-Y_{T'}^{\varphi_m}\right|
^{2}+2\int_{t}^{T'}\big( Y_{s}^{\varphi_n}-Y_{s}^{\varphi_m}\big)
\big(\varphi_{n}(s,Y_{s}^{\varphi_n},Z_{s}^{\varphi_n})-
\varphi_{m}(s,Y_{s}^{\varphi_m},Z_{s}^{\varphi_m})\big)
ds
\\& -2\int_{t}^{T'}\langle Y_{s}^{\varphi_n}-Y_{s}^{\varphi_m},
\quad \left(
Z_{s}^{\varphi_n}-Z_{s}^{\varphi_m}\right)dW_{s}\rangle.
\end{align*}
\\
For \;$N\in \N^\star$ we set, $ \Delta_{t}:=\left|
Y_{t}^{\varphi_n}-Y_{t}^{\varphi_m}\right| ^{2}+ (A_N)^{-1}. $
\\
Let \;$C>0$ and $1<\beta < \min
\{(3-\frac{2}{\overline\alpha}),2\}$. It\^o's formula shows that,
\begin{align*}
&e^{Ct}\Delta_{t}^{\beta\over 2}
+C\int_{t}^{T'}e^{Cs}\Delta_{s}^{\frac{\beta}{2}}ds\\& =
e^{CT'}\Delta_{T'}^{\beta\over 2}
+\beta\int_{t}^{T'}e^{Cs}\Delta_{s}^{\frac{\beta}{2}-1} \big(
Y_{s}^{\varphi_n}-Y_{s}^{\varphi_m}\big)
\big(\varphi_{n}(s,Y_{s}^{\varphi_n},Z_{s}^{\varphi_n})-
\varphi_{m}(s,Y_{s}^{\varphi_m},Z_{s}^{\varphi_m})\big)
ds
\\ & -\frac{\beta}{2}\int_{t}^{T'}e^{Cs}\Delta_{s}^{\frac{\beta}{2}-1}\left|
Z_{s}^{\varphi_n}-Z_{s}^{\varphi_m}\right|^{2}ds
-\beta\int_{t}^{T'}e^{Cs}\Delta_{s}^{\frac{\beta}{2}-1}\langle
Y_{s}^{\varphi_n}-Y_{s}^{\varphi_m}, \quad \left(
Z_{s}^{\varphi_n}-Z_{s}^{\varphi_m}\right)dW_{s}\rangle
\\ &-\beta(\frac{\beta}{2}-1)\int_{t}^{T'}e^{Cs}
\Delta_{s}^{\frac{\beta}{2}-2}
\left((Y_{s}^{\varphi_n}-Y_{s}^{\varphi_m})(Z_{s}^{\varphi_n}-
Z_{s}^{\varphi_m})\right)^{2}ds
\end{align*}
Put $\Phi(s)=|Y_{s}^{\varphi_n}| + |Y_{s}^{\varphi_m}|+
|Z_{s}^{\varphi_n}| + |Z_{s}^{\varphi_m}|$. Then
\begin{align*}
&
e^{Ct}\Delta_{t}^{\beta\over 2}
+C\int_{t}^{T'}e^{Cs}\Delta_{s}^{\frac{\beta}{2}}ds
\\& =
e^{CT'}\Delta_{T'}^{\beta\over 2}
-\beta\int_{t}^{T'}e^{Cs}\Delta_{s}^{\frac{\beta}{2}-1}\langle
Y_{s}^{\varphi_n}-Y_{s}^{\varphi_m}, \quad \left(
Z_{s}^{\varphi_n}-Z_{s}^{\varphi_m}\right)dW_{s}\rangle
\\ &
-\frac{\beta}{2}\int_{t}^{T'}e^{Cs}\Delta_{s}^{\frac{\beta}{2}-1}\left|
Z_{s}^{\varphi_n}-Z_{s}^{\varphi_m}\right|^{2}ds
\\ &
+\beta\frac{(2-\beta)}{2}\int_{t}^{T'}e^{Cs}
\Delta_{s}^{\frac{\beta}{2}-2}
\left((Y_{s}^{\varphi_n}-Y_{s}^{\varphi_m})(Z_{s}^{\varphi_n}-
Z_{s}^{\varphi_m})\right)^{2}ds
\\ &
+J_{1}+J_{2}+J_{3}+J_{4},
\end{align*}
where \begin{align*}
 & J_{1}:=\beta\int_{t}^{T'}e^{Cs}\Delta_{s}^{\frac{\beta}{2}-1}
\big( Y_{s}^{\varphi_n}-Y_{s}^{\varphi_m}\big)
\big(\varphi_{n}(s,Y_{s}^{\varphi_n},Z_{s}^{\varphi_n})-
\varphi_{m}(s,Y_{s}^{\varphi_m},Z_{s}^{\varphi_m})\big)
\1_{\{\Phi(s)>N\}} ds.
 \\ & J_{2}:=\beta\int_{t}^{T'}e^{Cs}\Delta_{s}^{\frac{\beta}{2}-1}
\big( Y_{s}^{\varphi_n}-Y_{s}^{\varphi_m}\big)
\big(\varphi_{n}(s,Y_{s}^{\varphi_n},Z_{s}^{\varphi_n})-
\varphi(s,Y_{s}^{\varphi_n},Z_{s}^{\varphi_n})\big)
\1_{\{\Phi(s)\leq N\}} ds.
\\ & J_{3}:=\beta\int_{t}^{T'}e^{Cs}\Delta_{s}^{\frac{\beta}{2}-1}
\big( Y_{s}^{\varphi_n}-Y_{s}^{\varphi_m}\big)
\big(\varphi(s,Y_{s}^{\varphi_n},Z_{s}^{\varphi_n})-
\varphi(s,Y_{s}^{\varphi_m},Z_{s}^{\varphi_m})\big)
\1_{\{\Phi(s)\leq N\}}ds.
\\ & J_{4}:=\beta\int_{t}^{T'}e^{Cs}\Delta_{s}^{\frac{\beta}{2}-1}
\big( Y_{s}^{\varphi_n}-Y_{s}^{\varphi_m}\big)
\big(\varphi(s,Y_{s}^{\varphi_m},Z_{s}^{\varphi_m})-
\varphi_{m}(s,Y_{s}^{\varphi_m},Z_{s}^{\varphi_m})\big)
\1_{\{\Phi(s)\leq N\}} ds.
\end{align*}
We shall estimate $J_{1}$, $J_{2}$, $J_{3}$, $J_{4}$. Let $\kappa
=3-\frac{2}{\overline\alpha}-\beta$. Since $\frac
{(\beta-1)}{2}+\frac{\kappa}{2}+\frac{1}{\overline\alpha}=1$, we
use H\" older inequality to obtain
\begin{align*}
J_{1} &\leq \beta e^{CT'} \dfrac{1}{N^\kappa}\int_{t}^{T'}
\Delta_{s}^{\frac{\beta-1}{2}}{\Phi^\kappa(s)}
|\varphi_{n}(s,Y_{s}^{\varphi_n},Z_{s}^{\varphi_n})-\varphi_{m}(s,Y_{s}^{\varphi_m},Z_{s}^{\varphi_m}|ds
\\ &\leq
\beta e^{CT'} \dfrac{1}{N^\kappa} \left[\int_{t}^{T'} \Delta_{s}
ds \right]^{\frac{\beta-1}{2}} \left[\int_{t}^{T'}{\Phi(s)}^2 ds
\right]^{\frac{\kappa}{2}}
\\ &
\times\left[\int_{t}^{T'}|\varphi_{n}(s,Y_{s}^{\varphi_n},Z_{s}^{\varphi_n})-
\varphi_{m}(s,Y_{s}^{\varphi_m},Z_{s}^{\varphi_m}|^{\overline\alpha}ds
\right]^{\frac{1}{\overline\alpha}}.
\end{align*}
Since $|Y_{s}^{\varphi_n}-Y_{s}^{\varphi_m}|\leq
\Delta_s^\frac{1}{2}$, it easy to see that
\begin{align*}
J_{2}+ J_{4} & \leq 2\beta e^{CT'} [ 2N^2+\nu_1]^{\frac{\beta-1}{2}}
\bigg[\int_{t}^{T'}\sup_{|y|,|z|\leq
N}|\varphi_{n}(s,y,z)-\varphi(s,y,z)|ds
\\  &
+\int_{t}^{T'}\sup_{|y|,|z|\leq
N}|\varphi_{m}(s,y,z)-\varphi(s,y,z)|ds\bigg].
\end{align*}
Using assumption \bf{(H3)}, we get
\begin{align*}
J_{3} &\leq \beta M_2
\int_{t}^{T'}e^{Cs}\Delta_{s}^{\frac{\beta}{2}-1}
 \bigg[|Y_{s}^{\varphi_n}-Y_{s}^{\varphi_m}|^{2}\log
A_{N}
\\ &
+ \frac{\log A_N}{A_N} +
|Y_{s}^{\varphi_n}-Y_{s}^{\varphi_m}||Z_{s}^{\varphi_n}-Z_{s}^{\varphi_m}|\sqrt{\log
A_{N}} \bigg]\1_{\{\Phi(s)<N\}}ds
\\ & \leq
\beta M_2\int_{t}^{T'}e^{Cs}\Delta_{s}^{\frac{\beta}{2}-1}
 \bigg[\Delta_{s}\log A_{N}+|Y_{s}^{\varphi_n}-Y_{s}^{\varphi_m}||Z_{s}^{\varphi_n}-Z_{s}^{\varphi_m}|\sqrt{\log
A_{N}}  \bigg]\1_{\{\Phi(s)\leq N\}}ds.
\end{align*}
%%%%%%%%%%%%%%%%%%%%%%%%%%%%%%%%%%%%%%%%%%%%%%%%%%%%%%%%%%%%%%%%%%%%%%%%%%
We choose $C=C_N=\dfrac{2M_2^2\beta}{\beta -1} \log A_{N}$, then
we use Lemma \ref{lem3} to show that
 $$\begin{array}{l}
e^{C_Nt}\Delta_{t}^{\beta\over 2} + \dfrac{\beta(\beta-1)}{4}
\int_{t}^{T'}e^{C_Ns}\Delta_{s}^{\frac{\beta}{2}-1}\left|
Z_{s}^{\varphi_n}-Z_{s}^{\varphi_m}\right|^{2}ds
\\ \leq e^{C_NT'}\Delta_{T'}^{\beta\over 2} -\beta\int_{t}^{T'}
e^{C_Ns}\Delta_{s}^{\frac{\beta}{2}-1}\langle
Y_{s}^{\varphi_n}-Y_{s}^{\varphi_m}, \quad \left(
Z_{s}^{\varphi_n}-Z_{s}^{\varphi_m}\right)dW_{s}\rangle
\\ +\beta e^{C_NT'} \dfrac{1}{N^\kappa}
\left[\int_{t}^{T'} \Delta_{s} ds \right]^{\frac{\beta-1}{2}}
\times\left[\int_{t}^{T'}{\Phi(s)}^2 ds \right]^{\frac{\kappa}{2}}
\\
\times\left[\int_{t}^{T'}|\varphi_{n}(s,Y_{s}^{\varphi_n},Z_{s}^{\varphi_n})-\varphi_{m}(s,Y_{s}^{\varphi_m},Z_{s}^{\varphi_m}|^{\overline\alpha}
\1_{\{\Phi(s)>N\}} ds \right]^{\frac{1}{\overline\alpha}}
\\  +\beta e^{C_NT'} [2N^2+\nu_1]^{\frac{\beta-1}{2}}
\bigg[\int_{t}^{T'}\sup_{|y|,|z|\leq
N}|\varphi_{n}(s,y,z)-\varphi(s,y,z)|ds
\\ +\int_{t}^{T'}\sup_{|y|,|z|\leq
N}|\varphi_{m}(s,y,z)-\varphi(s,y,z)|ds\bigg]
\end{array}$$
Burkholder's inequality and H\" older's inequality (since $\frac
{(\beta-1)}{2}+\frac{\kappa}{2}+\frac{1}{\overline\alpha}=1$)
allow us to show that there exists a universal constant $\ell>0$
such that $\forall \delta'>0$,
\begin{align*}
& \E \sup_{(T'-\delta')^{+}\leq t \leq T'}\left[
e^{C_Nt}\Delta_{t}^{\beta\over 2}\right] + \E
\int_{(T'-\delta')^{+}}^{T'}e^{C_Ns}\Delta_{s}^{\frac{\beta}{2}-1}\left|
Z_{s}^{\varphi_n}-Z_{s}^{\varphi_m}\right|^{2}ds
\\&\leq \frac{\ell}{\beta -1} e^{C_NT'}\bigg\{ \E \left[\Delta_{T'}^{\beta\over 2}\right]
+\dfrac{\beta}{N^\kappa} \left[\E \int_{0}^{T} \Delta_{s} ds
\right]^{\frac{\beta-1}{2}} \left[\E \int_{0}^{T}{\Phi(s)^2} ds
\right]^{\frac{\kappa}{2}}
\\ & \times\bigg[\E
\int_{0}^{T}|\varphi_{n}(s,Y_{s}^{\varphi_n},Z_{s}^{\varphi_n})-\varphi_{m}(s,Y_{s}^{\varphi_m},Z_{s}^{\varphi_m}|^{\overline\alpha}
ds \bigg]^{\frac{1}{\overline\alpha}}
\\ & +\beta [2N^2+\nu_1]^{\frac{\beta-1}{2}}
\E \bigg[\int_{0}^{T}\sup_{|y|,|z|\leq
N}|\varphi_{n}(s,y,z)-\varphi(s,y,z)|ds
\\ &+\int_{0}^{T}\sup_{|y|,|z|\leq
N}|\varphi_{m}(s,y,z)-\varphi(s,y,z)|ds\bigg] \bigg\}.
\end{align*}
%%%%%%%%%%%%%%%%%%%%%%%%%%%%%%%%%%%%%%%%%%%%%%%%%%%%%%%%%%%%%%%%%%%%%%%%%%%%%%%%%%%%%%
We use
%Lemma \ref{estimateY}, Lemma \ref{estimateZ}, Lemma  \ref{estimatevarphi},
Lemma \ref{lem1} and
Lemma} \ref{lem2} to obtain, $\forall N>R$,
\begin{align*}
& \E \sup_{(T'-\delta')^{+}\leq t \leq T'}\vert
Y_{t}^{\varphi_n}-Y_{t}^{\varphi_m}\vert^\beta + \E
\int_{(T'-\delta')^{+}}^{T'}\dfrac{\left|
Z_{s}^{\varphi_n}-Z_{s}^{\varphi_m}\right|^{2}}{\left(\vert
Y_{s}^{\varphi_n}-Y_{s}^{\varphi_m}\vert^{2}+ \nu_{R}
\right)^{\frac{2-\beta}{2}}}ds
\\&\leq \frac{\ell}{\beta -1} e^{C_N\delta'} \bigg\{(A_N)^{-\beta\over 2}
+\beta \dfrac{2
K_{3}^{\frac{1}{\overline\alpha}}}{N^{\kappa}}\left(4TK_2
+T\ell\right)^{\frac{\beta-1}{2}}\left(8TK_2
+8K_1\right)^{\frac{\kappa}{2}}
\\ & + \E \vert
Y_{T'}^{\varphi_n}-Y_{T'}^{\varphi_m}\vert^\beta +\beta [2N^2+\nu_1
]^{\frac{\beta-1}{2}} \big[\rho_{N}(\varphi_n -
\varphi)+\rho_{N}(\varphi_m - \varphi)\big]\bigg\}
\\ & \leq
\frac{\ell}{\beta -1} e^{C_N\delta'}\E \vert
Y_{T'}^{\varphi_n}-Y_{T'}^{\varphi_m}\vert^\beta+ \frac{\ell}{\beta
-1}
\dfrac{A_N^{\frac{2M_2^2\delta'\beta}{\beta-1}}}{(A_N)^{\frac{\beta}{2}}}
\\
&+ \frac{2\ell}{\beta -1}\beta
K_{3}^{\frac{1}{\overline\alpha}}\left(4TK_2
+T\ell\right)^{\frac{\beta-1}{2}}\left(8TK_2
+8K_1\right)^{\frac{\kappa}{2}}\dfrac{A_N^{\frac{2M_2^2\delta'\beta}{\beta-1}}}
{(A_N)^{\frac{\kappa}{r}}}
\\ & +\frac{2\ell}{\beta -1} e^{C_N\delta'}\beta [2N^2+\nu_1
]^{\frac{\beta-1}{2}} \big[\rho_{N}(\varphi_n -
\varphi)+\rho_{N}(\varphi_m - \varphi)\big].
\end{align*}
%%%%%%%%%%%%%%%%%%%%%%%%%%%%%%%%%%%%%%%%%%%%%%%%%%%%%%%%%%%%%%%%%%%%%%%%%%%%%%%%%%%%
Hence for $ \delta' < (\beta
-1)\min\left(\frac{1}{4M_2^2},\frac{\kappa}{2rM_2^2\beta}\right)$
we derive
\begin{equation*}
\dfrac{A_N^{\frac{2M_2^2\delta'\beta}{\beta-1}}}{(A_N)^{\frac{\beta}{2}}}\longrightarrow_{N\rightarrow\infty}
0
\end{equation*}
 and
\begin{align*}
  \displaystyle\dfrac{A_N^{\frac{2M_2^2\delta'\beta}{\beta-1}}}
{(A_N)^{\frac{\kappa}{r}}}
\displaystyle\longrightarrow_{N\rightarrow\infty} 0.
\end{align*}
Passing to the limits first on $n$ and next on $N$, and using
assertion $(c)$ of lemma \ref{lem1}. \eop

\begin{remark} To deal with the case which take account of the
process $v_t$ appearing in assumption \bf{(H3)}, it suffices to take \ $\Phi(s):=|Y_{s}^{1}| + |Y_{s}^{2}|+
|Z_{s}^{1}| + |Z_{s}^{2}|+v_s$ \ in the proof of Lemma \ref{lem4}.
\end{remark}
%%%%%%%%%%%%%%%%%%%%%%%%%%%%%%%%%%%%%%%%%%%%%%%%%%%%%%%%%%%%%%%%%%%%%%%%%%%%%%%%%%%%%%
%                   Proof of theorem 2.1
%%%%%%%%%%%%%%%%%%%%%%%%%%%%%%%%%%%%%%%%%%%%%%%%%%%%%%%%%%%%%%%%%%%%%%%%%%%%%%%%%%%%
\bop \bf {of Theorem \ref{unique}} Taking successively $T'=T$,
$T'=(T-\delta')^+$, $T'=(T-2\delta')^{+}...$ in Lemma \ref{lem4}, we
obtain, for every $\beta\in ]1,\quad
\min\left(3-\dfrac{2}{\overline\alpha}, 2\right)[$
\begin{align*}
\lim_{n,m\rightarrow +\infty}\left( \E \sup_{0\leq t \leq T}\vert
Y_{t}^{\varphi_n}-Y_{t}^{\varphi_m}\vert^\beta  + \E
\int_{0}^{T}\dfrac{\left|
Z_{s}^{\varphi_n}-Z_{s}^{\varphi_m}\right|^{2}}{\left(\vert
Y_{s}^{\varphi_n}-Y_{s}^{\varphi_m}\vert^{2}+\nu_{R}\right)^{\frac{2-\beta}{2}}}ds\right)=
0.
\end{align*}
But by Schwarz inequality we have \vskip 0.2cm\noindent
$$
\E\dint_{0}^{T}\vert Z_{s}^{\varphi_n} - Z_{s}^{\varphi_m}\vert ds
\leq \bigg(\E \dint_{0}^{T}\dfrac{\left|
Z_{s}^{\varphi_n}-Z_{s}^{\varphi_m}\right|^{2}}{\left(\vert
Y_{s}^{\varphi_n}-Y_{s}^{\varphi_m}\vert^{2}+\nu_{R}\right)^{\frac{2-\beta}{2}}}ds\bigg)^\frac{1}{2}
\bigg(\E \dint_{0}^{T}{\left(\vert
Y_{s}^{\varphi_n}-Y_{s}^{\varphi_m}\vert^{2}+\nu_{R}\right)^{\frac{2-\beta}{2}}}ds\bigg)^\frac{1}{2}
$$
\par\noindent
Since $\beta >1$, Lemma \ref{lem2} allows us to show that
\begin{align*}
\lim_{n\rightarrow +\infty} \left(\E\sup_{0\leq t\leq T}\vert
Y_{t}^{\varphi_n} - Y_{t}\vert^{\beta}+\E\dint_{0}^{T}\vert
Z_{s}^{\varphi_n} - Z_{s}\vert ds\right) = 0.
\end{align*}
In particular, there exists a subsequence, which we still denote
$(Y^{\varphi_n}$, $Z^{\varphi_n})$, such that
\begin{align*}
\lim_{n\rightarrow +\infty}\left(\vert Y_{t}^{\varphi_n} -
Y_{t}\vert+\vert Z_{t}^{\varphi_n} - Z_{t}\vert\right) = 0\quad
a.e.\; (t,\omega).
\end{align*}
On the other hand
\begin{align*}
&\E\dint_{0}^{T}\vert \varphi_{n}(s,Y_{s}^{\varphi_n},
Z_{s}^{\varphi_n}) -\varphi(s,Y_{s}^{\varphi_n},
Z_{s}^{\varphi_n})\vert ds
\\ & \leq \E\dint_{0}^{T}\vert \varphi_{n}(s,Y_{s}^{\varphi_n}, Z_{s}^{\varphi_n})-f(s,Y_{s}^{\varphi_n},
Z_{s}^{\varphi_n})\vert \1_{\{\vert Y_{s}^{\varphi_n}\vert +\vert
Z_{s}^{\varphi_n}\vert\leq N\}} ds
\\ & +\E\dint_{0}^{T}\vert \varphi_{n}(s,Y_{s}^{\varphi_n}, Z_{s}^{\varphi_n})-f(s,Y_{s}^{\varphi_n},
Z_{s}^{\varphi_n})\vert\dfrac{(\vert Y_{s}^{\varphi_n}\vert +\vert
Z_{s}^{\varphi_n}\vert)^{(2-\frac{2}{\overline\alpha})}}{N^{(2-\frac{2}{\overline\alpha})}}
\1_{\{\vert Y_{s}^{\varphi_n}\vert +\vert Z_{s}^{\varphi_n}\vert\geq
N\}} ds
\\ & \leq \rho_{N}(\varphi_n - \varphi)
+\dfrac{2K_{3}^{\frac{1}{\overline\alpha}}\left[TK_{2}
+K_1\right]^{1-\frac{1}{\overline\alpha}}}{N^{(2-\frac{2}{\overline\alpha})}}.
\end{align*}
Passing to the limit first on $n$ and next on $N$  we obtain
\begin{align*}
\lim_{n}E\dint_{0}^{T}\vert \varphi_{n}(s,Y_{s}^{\varphi_n},
Z_{s}^{\varphi_n}) -\varphi(s,Y_{s}^{\varphi_n},
Z_{s}^{\varphi_n})\vert ds =0.
\end{align*}
Finally, we use $\bf{(H.1)}$, Lemma \ref{lem1} and Lemma
\ref{lem2} to show that,
\begin{align*}
\lim_{n}E\dint_{0}^{T}\vert \varphi_{n}(s,Y_{s}^{\varphi_n},
Z_{s}^{\varphi_n}) -\varphi(s,Y_{s}, Z_{s})\vert ds =0.
\end{align*}
The existence is proved.

%%%%%%%%%%%%%%%%%%%%%%%%%%%%%%%%%%%%%%%%%%%%%%%%%%%%%%%%%%%%%%%%%%%%
%                         Uniqueness                               %
%%%%%%%%%%%%%%%%%%%%%%%%%%%%%%%%%%%%%%%%%%%%%%%%%%%%%%%%%%%%%%%%%%%%%%
\vskip 0.2cm\noindent
 \bf{Uniqueness.} Let $(Y,Z)$ and $(Y',Z')$ be two solutions  of
equation $(E^{f})$. Arguing as previously one can show that:
\\
for every $R>2$, $\beta \in
]1,\min\left(3-\dfrac{2}{\overline\alpha},2\right)[$, $\delta' <
(\beta
-1)\min\left(\frac{1}{4M_2^2},\frac{3-\frac{2}{\overline\alpha}-\beta}{2rM_2^2\beta}\right)$
and $\varepsilon>0$\\
 there exists $N_{0}>R$ such that  for all $N>N_{0}$, $\forall T'\leq T$
\begin{align*}
& \E \sup_{(T'-\delta')^{+}\leq t \leq T'}\vert
Y_{t}-Y_{t}^{'}\vert^\beta  + \E
\int_{(T'-\delta')^{+}}^{T'}\dfrac{\left|
Z_{s}-Z_{s}^{'}\right|^{2}}{\left(\vert
Y_{s}-Y_{s}^{'}\vert^{2}+\nu_{R}\right)^{\frac{2-\beta}{2}}}ds
\\ &\qquad\qquad\qquad\leq
\varepsilon +\frac{\ell}{\beta -1} e^{C_N\delta'}\E \vert
Y_{T'}-Y_{T'}^{'}\vert^\beta.
\end{align*}
Again, taking successively $T'=T$, $T'=(T-\delta')^+$,
$T'=(T-2\delta')^{+}...$, we establish the uniqueness of solution.
Theorem \ref{unique} is proved. \eop
%%%%%%%%%%%%%%%%%%%%%%%%%%%%%%%%%%%%%%%%%%%%%%%%%%%%%%%%%%%%%%%%%%%%
%                       Proof of stability Theorem 2.2
%%%%%%%%%%%%%%%%%%%%%%%%%%%%%%%%%%%%%%%%%%%%%%%%%%%%%%%%%%%%%%%%%%%%

\vskip 0.3cm \bop\bf{of Theorem \ref{stability}.} Also as in the proof
of Theorem \ref{unique}, we show that,
 \\
 For every $R>2$, $\beta \in ]1,\min\left(3-\dfrac{2}{\overline\alpha},2\right)[$,
 $\delta' <
 (\beta-1)\min\left(\frac{1}{4M_2^2},\frac{3-\frac{2}{\overline\alpha}-\beta}{2rM_2^2\beta}\right)$
 and
 $\varepsilon>0$, there exists $N_{0}>R$ such that for all $N>N_{0}$, for all $ T'\leq T$:
\begin{align*}
& \limsup_{n\rightarrow +\infty}\E \sup_{(T'-\delta')^{+}\leq t
\leq T'}\vert Y_{t}^n-Y_{t}\vert^\beta  + \E
\int_{(T'-\delta')^{+}}^{T'}\dfrac{\left|
Z_{s}^n-Z_{s}\right|^{2}}{\left(\vert
Y_{s}^n-Y_{s}\vert^{2}+\nu_{R}\right)^{\frac{2-\beta}{2}}}ds
\\ &\qquad\qquad\qquad\leq
\varepsilon +\frac{\ell}{\beta -1} e^{C_N\delta'}
\limsup_{n\rightarrow +\infty}\E \vert Y_{T'}^n-Y_{T'}\vert^\beta.
\end{align*}
Again as in the proof of Theorem \ref{unique}, taking successively
$T'=T$, $T'=(T-\delta')^+$, $T'=(T-2\delta')^{+}...$, we establish
the convergence in the whole interval $[0,T]$. In particular, we
have for every $q<2$, \quad
$
\lim_{n\rightarrow +\infty}\left(\vert Y^n - Y\vert^q\right)=0
$
\quad and \quad $\lim_{n\rightarrow +\infty}\left(\vert Z^n -
Z\vert^q\right) = 0 $ \quad \it{in measure} $P\times dt$. Since
$(Y^n)$ and $(Z^n)$ are square integrable, the proof is finished
by using an uniform integrability argument. Theorem \ref{stability} is
proved. \eop

\section{Application to stochastic and control}
\medskip
In all the following $\Omega={\cal C}([0,T],\R^m)$ is the space of
continuous
functions from $[0,T]$ to $\R^m$. %${\cal P}$ is the Wiener measure
%on $\Omega$, $B_t(w)=w(t)$ for all $(w,t)\in \Omega \times[0,T]$,
%and ${\cal F}_{t}^{0}:=\sigma \{B_{s},s\leq t\}$. The process
%$(B_t,t\leq T)$ is an $({\cal F}_{t}^{0},{\cal P})$ Brownian motion.
%we shall denote by $\mathbb{F}=\{{\cal F}_t\}_{t\leq T}$ the ${\cal
%P}$-augmentation of this natural
%filtration.\\

Let us consider a mapping $\sigma:\, (t,w)\in[0,T]\times \Omega
\rightarrow \sigma(t,w)\in\R^{m}\bigotimes \R^{m}$ satisfying the
following:

\medskip
$(1.1)$ $\sigma$ is {\cal P}-measurable.

$(1.2)$ There exists a constant $C$ such that
$|\sigma(t,w)-\sigma(t,w')|\leq C||w-w'||_t$ and $|\sigma(t,w)|\leq
C(1+||w||_t)$, where for any $w,w'\in\Omega^2$ and $t\leq T,\,
||w||_t=\sup\limits_{s\leq t}|w_s|$.

$(1.3)$ For any $(t,w)\in [0,T]\times \Omega$, the matrix
$\sigma(t,w)$ is invertible and $|\sigma^{-1}(t,w)|\leq
C$ for some constants $C$.\\
\medskip

Let $x_0\in\R^m$ and $x=(x_t)_{t\leq T}$ be the solution of the
following standard functional differential equation:
\begin{equation}\label{EDS}
x_t=x_0+\int_0^t\sigma(s,x)dB_s,\quad t\leq T;
\end{equation}
the process $(x_t)_{t\leq T}$ exists, since $\sigma$ satisfies
$(1.1)-(1.3)$ (see,e,g.,\cite{RY} page 375. Moreover,
\begin{equation}\label{estim-eds}
\E[(||x||_T)^n]<+\infty,\quad \forall n\in[1,+\infty[ (\cite{KS},pp.
306).
\end{equation}
\subsection{Stochastic control of diffusions}

Let $A$ be a compact metric space and ${\cal U}$ be the space of
${\cal P}$-measurable processes $u:=(u_t)_{t\leq T}$ with value in
$A$. Let $f:[0,T]\times \Omega \times A\rightarrow \R^m$ be such
that:

\medskip
$(1.4)$ For each $a\in A$, the function $(t,w)\rightarrow f(t,w,a)$
is predictable.

$(1.5)$ For each $(t,w)$, the mapping $a\rightarrow f(t,w,a)$ is
continuous.

 $(1.6)$ There exists a real constant $K>0$ such that
\begin{equation}\label{generateur-cond}
|f(t,w,a)|\leq K(1+||w||_t),\quad \forall 0\leq t \leq T,\,w\in
\Omega,\,a\in A.
\end{equation}
\medskip

For any given admissible control strategy $u\in {\cal U}$, the
exponential process

$$\Lambda^u_t=exp\{\int_0^T\sigma^{-1}(s,x)f(s,x,u_s)dB_s-\frac{1}{2}\int_0^T|\sigma^{-1}(s,x)f(s,x,u_s)|^2ds\}$$
$0\leq t \leq T$, is a martingale under all these assumptions;
namely, $\E[\Lambda^u_T]=1$ (see Karatzas and Shreve (1991), pages
191 and 200 for this result). Then the Girsanov theorem guarantees
that the process
\begin{equation}\label{chang-prob}
B^u_t=B_t-\int_0^t\sigma^{-1}(s,x)f(s,x,u_s)ds,\quad 0\leq t \leq T,
\end{equation}
is a Brownian motion with respect to the filtration ${\cal F}_t$,
under the new probability measure
$$P^u(B)=\E[\Lambda^u_T.\ind_{B}],\quad B\in {\cal F}_T,$$
which is equivalent to $P$. It is now clear from the equations
(\ref{EDS}) and (\ref{chang-prob}) that
\begin{equation}\label{eds-contr}
x_t=x_0+\int_0^tf(s,x,u_s)ds+\int_0^t\sigma(s,x)dB^u_s,\quad 0\leq t
\leq T,
\end{equation}
holds almost surely. This will be our model for a controlled
stochastic functional differential equation, with the control
appearing only in the drift term.

In order to specify the objective of our stochastic game of control
and stopping. Let us now consider the followings:
\medskip

\indent $(1.6)$ $h:[0,T]\times \Omega \times A\rightarrow \R$ is
measurable and for each $(t,w)$ the mapping $a\rightarrow h(t,w,a)$
is continuous. In addition there exists a real constant $K>0$ such
that
\begin{equation}\label{hgenerateur}
|h(t,w,a)|\leq K(1+||w||_t),\quad \forall 0 \leq t \leq T,\,w\in
\Omega,\,a\in A.
\end{equation}

$(1.7)$ %$Q=(Q_t)_{t\leq T}$ is a process of ${\cal S}^p$,
$g_1:[0,T]\times \Omega \rightarrow \R$ and is continuous function
and there exists a real positive constant $C$  such that:
\begin{equation} \label{polycond} |g_1(t,w)|\leq
C(1+||w||_t),\, \, \forall (t,w)\in [0,T]\times \Omega .
\end{equation}
% Moreover we assume that $Q_t\geq g_1(t,x)$ and $Q_T=g_2(T,x)\geq
% g_1(T,x)$, $\forall (t,x)\in [0,T]\times \R^m$.
\medskip

We shall study a stochastic control with one player. The controller,
who chooses an admissible control strategy $u\in {\cal U}$ to
minimize this  amount
\begin{equation}\label{amount}\int_0^{T}h(s,x,u_s)ds+g_1(T,x_T).\end{equation} It is thus in the best interest of the
controller  to make the amount (\ref{amount}) as small
 as possible, at least on the average. We
are thus led to a stochastic control, with
\begin{equation}\label{jeux}J(u)=\E^u[\int_0^{T}h(s,x,u_s)ds+g_1(T,x_T)].\end{equation} The problem we are interested in is finding
an intervention strategies $u^*$, for controller such that for any
$u\in {\cal U}$, we have
$$J(u^*)\leq J(u).$$ Then $u^*$ is called an optimal control for the problem. Now let us set
\begin{equation}\label{Hamiltonien}H(t,x,z,u_t)=z\sigma^{-1}(t,x)f(t,x,u_t)+h(t,x,u_t)\quad \forall(t,x,z,u_t)\in [0,T]\times\Omega\times\R^{m}\times
A.\end{equation} The function $H$ is called the Hamiltonian
associated with stochastic control such that:
\medskip

$(2.1)$ $ \forall z\in \R^{m}$, the process $(H(t,x,z,u_t))_{t\leq
T}$ is ${\cal P}$-measurable.

%                            lemma monotonie
%%%%%%%%%%%%%%%%%%%%%%%%%%%%%%%%%%%%%%%%%%%%%%%%%%%%%%%%%%%

\begin{lemma}\label{monotonie}
The Hamiltonian $H$ satisfies  \bf{(H.2)} and \bf{(H.3)}.
\end{lemma}

\vskip 0.2cm

%%%%%%%%%%%%%%%%%%%%%%%%%%%%%%%%%%%%%%%%%%%%%%%%%%%%%%%%%%%
%                            proof lemma monotonie
%%%%%%%%%%%%%%%%%%%%%%%%%%%%%%%%%%%%%%%%%%%%%%%%%%%%%%%%%%%

\bop For \bf{(H.2)}, it is not difficult to show that for every \
$(t,x,z,u_t)\in [0,T]\times\Omega\times\R^{m}\times A$ and $\mid
z\mid$ large enough, there exist a constants $C$ and $c_0$ such
that:
 \begin{equation}\label{majora_gener}|H(t,x,z,u_t)|\leq
Cexp(||x||_t)+c_0|z|\ln^{\frac{1}{2}}(|z|).\end{equation}
\medskip

To prove that $H$ satisfies assumption \bf{(H.3)}, it is enough to
take $v_t := \exp\vert f(t,x,u_t)\vert$. Indeed, we have

\begin{align*}
\big( y-y^{\prime}\big) \big(H(t,x,z,u_t)-H(t,x,z',u_t)\big) \1_{\{
e^{\vert f(t,x,u_t)\vert^2}\leq N\}} & \leq \mid
y-y^{\prime}\mid\mid z-z^{\prime}\vert
f(t,x,u_t)\vert\1_{\{ \vert f(t,x,u_t)\vert^2\leq \log N\}} \\
& \leq \mid y-y^{\prime}\mid\mid z-z^{\prime}\vert \sqrt{\log A_{N}}
\end{align*}
To complete the proof, we shall show that $ \exp\vert
f(t,x,u_t)\vert^2$ belongs to $ L^{q}(\Omega\times [0, T]; \R_
 +))$ for some $q>0$. We have,

\begin{align*}
\E\int_0^T \exp(q\vert f(s,x,u_s)\vert^2)ds &\leq \E\int_0^T
\exp(2qK^2(1 + \sup_{s\leq T}\vert x_s\vert^2)ds \\
&\leq \exp(2qK^2T) \E\int_0^T \exp(2qK^2 \sup_{s\leq T}\vert x_s\vert^2)ds \\
&\leq  T\exp(2qK^2T)\E\exp(2qK^2 \sup_{s\leq T}\vert x_s\vert^2) \\
\end{align*}
And, since $\sigma$ is with linear growth, it is well known that
$E\exp(2qK^2 \sup_{s\leq T}\vert x_s\vert^2) < \infty$ for $q$ small
enough. \eop

\medskip

%%%%%%%%%%%%%%%%%%%%%%%%%%%%%%%%%%%%%%%%%%%%%%%%%%%%%%%%%%%
%                            lemma estimateY
%%%%%%%%%%%%%%%%%%%%%%%%%%%%%%%%%%%%%%%%%%%%%%%%%%%%%%%%%%%

To begin with let us define the notion of solution of the reflected
BSDE associated with the triple  $(H,g_2,g_1)$ which we consider
throughout this paper.

In order to construct a stochastic control, we need to do is find an
admissible control strategy $u^*()\in {\cal U}$ for our stochastic control.\\
The Hamiltonian function defined in (\ref{Hamiltonien}) attains its
infimum over the set $A$ at some $u^* \equiv u^*(t,x,p)\in A$, for
any given $(t,x,p)\in[0,T]\times \Omega \times \R^m$, namely,
\begin{equation}\label{BN}
\inf\limits_{u\in A}H(t,x,u,p)=H(t,x,u^*(t,x,p),p).
\end{equation}
(This is the case, for instance, if the set $A$ is compact and the
mapping $u\rightarrow H(t,x,u,p)$ continuous.) Then it can be shown
(see Lemma 1 in Benes (1970), that the mapping $u^*:([0; T]\times
\Omega \times \R^m \rightarrow A$ can be selected to be ${\cal P}
\otimes {\cal
B}(\R^m)$-measurable.\\

Now let $H^*(t,x,z)=\inf\limits_{u\in A}H(t,x,u,z)$ where $x$ is the
solution of (\ref{EDS}). Let $(Y_t)_{t\leq T}$ be the process
constructed as in Theorem \ref{unique} with $(H^*,g_1)$. Using once
again Theorem \ref{unique}, there exists a unique pair $(Y_t,Z_t)_{t
\leq T}$ such that
\begin{equation}
\left\{
  \begin{array}{ll}
  (Y,Z)\in (\E, \vert\vert.\vert\vert);\\
\displaystyle Y_t = g_1(T,x_T)+ \int_t^{T} H^*(s,x,Z_s) ds
 - \int_t^{T} Z_s dB_s$, $t\in[0,T].\label{def1-rbsde}
  \end{array}
\right.
\end{equation}

 We are ready to give the main
result of this section.
\begin{theorem}\label{games}The admissible control $u^*$ is optimal for the stochastic control; i.e., it
satisfies $$J(u^*)=Y_0\leq J(u)\,\, \forall u\in{\cal U}.$$
Additionally, $Y_0$ is the value of the stochastic control, i.e.,
$$Y_0=\inf\limits_{u\in{\cal U}}J(u).
$$
\end{theorem}
\bop: Let us show that $Y_0=J(u^*)$. It follows that $$
\begin{array}{l} \displaystyle Y_0=g_1(T,x_T)+ \int_0^{T} H^*(s,x,Z_s) ds - \int_0^{T} Z_s dB_s\\
\quad\, \displaystyle  =g_1(T,x_T)+ \int_0^{T} h(s,x,u^*(s,x,Z_s))ds
- \int_0^{T} Z_s dB^{u^*}_s.
\end{array}$$
As $(\int_0^{t} Z_s dB_s)_{t\leq T}$ is an
$(\mathbb{F}_t,P^{u^*})$-martingale, taking expectation we get
$$\displaystyle Y_0=\E^{u^*}[Y_0]=\E^{u^*}[g_1(T,x_T)+ \int_0^{T} h(s,x,u^*(s,x,Z_s))
ds],$$ because $Y_0$ is $\mathbb{F}_0$-measurable, and hence
deterministic. Now $P \mbox{-a.s.}$, and also $P^{u^*} \mbox{-a.s.}$
(since they are equivalent probabilities). Then

$$
\begin{array}{l} \displaystyle
Y_0=J(u^*).\end{array}$$

 Next let $u\in{\cal U}$. Let us show that $Y_0\leq J(u)$.
 $$
\begin{array}{lll} \displaystyle Y_0=g_1(T,x_T)+ \int_0^T H^*(s,x,Z_s) ds - \int_0^{T} Z_s dB_s
\\
\quad\, \displaystyle \leq g_1(T,x_T)+ \int_0^T H(s,x,u,Z_s) ds -
\int_0^{T} Z_s dB_s
\\
\quad\, \displaystyle  =g_1(T,x_T)+\int_0^{T} h(s,x,u(s,x,Z_s))ds -
\int_0^{T} Z_s dB^{u}_s,
\end{array}$$
 Once more $(\int_0^{t} Z_s dB_s)_{t\leq
T}$ is an $(\mathbb{F}_t,P^{u})$-martingale; then taking the
expectation with respect to $P^{u}$ and taking into account the fact
that $Y_0$ is deterministic, we obtain
$$Y_0=\E^{u}[Y_0]\leq \E^{u}[g_1(T,x_T)+ \int_0^{T} h(s,x,u(s,x,Z_s)) ds],$$
then $Y_0\leq J(u)$. The proof is now complete. \eop

\subsection{Stochastic zero-sum differential games }

Let $A$ (resp. $B$) be a compact metric space and ${\cal U}$ (resp.
${\cal V}$) be the space of ${\cal P}$-measurable processes
$u:=(u_t)_{t\leq T}$ (resp. $v:=(v_t)_{t\leq T}$)with value in $A$
(resp. $B$). Let $f:[0,T]\times \Omega \times A \times B \rightarrow
\R^m$ be such that:

\medskip
$(1.4)$ For each $a\in A$ and $b\in B$, the function
$(t,x)\rightarrow f(t,w,a,b)$ is predictable.

$(1.5)$ For each $(t,w)$, the mapping $(a,b)\rightarrow f(t,w,a,b)$
is continuous.

 $(1.6)$ There exists a real constant $K>0$ such that
\begin{equation}\label{generateur-cond}
|f(t,w,a,b)|\leq K(1+||w||_t),\quad \forall 0\leq t \leq T,\,w\in
\Omega,\,a\in A, \,b\in B.
\end{equation}
\medskip

For any given admissible control strategy $(u,v)\in {\cal U}\times
{\cal V}$, the exponential process

$$\Lambda^{(u,v)}=exp\{\int_0^T\sigma^{-1}(s,x)f(s,x,u_s,v_s)dB_s-\frac{1}{2}\int_0^T|\sigma^{-1}(s,x)f(s,x,u_s,v_s)|^2ds\}$$
$0\leq t \leq T$, is a martingale under all these assumptions;
namely, $\E[\Lambda^{(u,v)}_T]=1$ (see Karatzas and Shreve (1991),
pages 191 and 200 for this result). Then the Girsanov theorem
guarantees that the process
\begin{equation}\label{chang-prob}
B^{(u_t,v_t)}=B_t-\int_0^t\sigma^{-1}(s,x)f(s,x,u_s,v_s)ds,\quad
0\leq t \leq T,
\end{equation}
is a Brownian motion with respect to the filtration ${\cal F}_t$,
under the new probability measure
$$P^{(u,v)}(B)=\E[\Lambda^{(u,v)}_T.\ind_{B}],\quad B\in {\cal F}_T,$$
which is equivalent to $P$. It is now clear from the equations
(\ref{EDS}) and (\ref{chang-prob}) that
\begin{equation}\label{eds-contr}
x_t=x_0+\int_0^tf(s,x,u_s,v_s)ds+\int_0^t\sigma(s,x)dB^{(u_s,v_s)},\quad
0\leq t \leq T,
\end{equation}
holds almost surely.

 Let us now consider the followings:
\medskip

\indent $(1.6)$ $h:[0,T]\times \Omega \times A \times B\rightarrow
\R$ is measurable and for each $(t,w)$ the mapping $(a,b)\rightarrow
h(t,w,a,b)$ is continuous. In addition there exists a real constant
$K>0$ such that
\begin{equation}\label{hgenerateur}
|h(t,w,a,b)|\leq K(1+||w||_t),\quad \forall 0 \leq t \leq T,\,w\in
\Omega,\,a\in A, \, b\in B.
\end{equation}

$(1.7)$ %$Q=(Q_t)_{t\leq T}$ is a process of ${\cal S}^p$,
$g_1:[0,T]\times\R^m\rightarrow \R$ and is continuous function and
there exists a real positive constant $C$  such that:
\begin{equation} \label{polycond} |g_1(t,w)|\leq
C(1+||w||_t),\, \, \forall (t,w)\in [0,T]\times \R^m. \end{equation}
% Moreover we assume that $Q_t\geq g_1(t,x)$ and $Q_T=g_2(T,x)\geq
% g_1(T,x)$, $\forall (t,x)\in [0,T]\times \R^m$.
\medskip

We shall study a stochastic zero-sum differential games. Then the
payoff corresponding to $u \in {\cal U}$ and $v \in {\cal V}$ is
\begin{equation}\label{jeux}J(u,v)=\E^{u,v}[\int_0^{T}h(s,x,u_s,v_s)ds+g_1(T,x_T)].\end{equation}
 where $u\in {\cal U}$ (resp. $v\in {\cal U}$) is the strategy of the
 first (resp. second) player. The first player looks for minimize $J(u,v)$, when the second
 looks for maximize the same $J(u,v)$.We are concerned by the problem of the existence
 of a saddle-point for this game, i.e the existence of an admissible contol
 $(u^*,v^*)$ which satisfies: $$J(u^*,v)\leq J(u^*,v^*)\leq J(u,v^*),\quad (u,v)\in {\cal U}\times {\cal V}$$ We introduce
 the hamiltonian function defined by:
\begin{equation}\label{Hamiltonien1}H(t,x,z,u_t,v_t)=z\sigma^{-1}(t,x)f(t,x,u_t,v_t)+h(t,x,u_t,v_t)\quad \forall(t,x,z,u_t,v_t)\in [0,T]\times\Omega\times\R^{m}\times
A\times B,\end{equation} and we suppose that the Isaacs' condition
is satisfied:
$$\mbox{(H)} \qquad\sup\limits_{v\in{\cal V}}\inf\limits_{u\in{\cal U}}H(t,x,p,u,v)=\inf\limits_{u\in{\cal U}}\sup\limits_{v\in{\cal V}}H(t,x,p,u,v)\quad \quad \forall(t,x,p)\in [0,T]\times\Omega\times\R^{m}$$
\medskip

$(2.1)$ $ \forall z\in \R^{m}$, the process
$(H(t,x,z,u_t,v_t))_{t\leq T}$ is ${\cal P}$-measurable.

Then using a selection theorem  \cite{B} we get easily the:
\begin{lemma}\label{BN}
$(H)$ is equivalent to the following assumption:

There exists $u^*(t,x,p), v^*(t,x,p)$ ${\cal P}\otimes {\cal
B}(\R^m)$- mesurable valued respectively in ${\cal U}$ and ${\cal
V}$ such that:
$$\begin{array}{l}H(t,x,p,u^*(t,x,p),v(t,x,p)) \\
\leq H(t,x,p,u^*(t,x,p),v^*(t,x,p)) \leq
H(t,x,p,u(t,x,p),v^*(t,x,p))\quad \forall u,v,t,x,p.\end{array}$$
Moreover $u^*$ and $v^*$ satisfy

$$H(t,x,p,u^*(t,x,p),v^*(t,x,p))=\sup\limits_{v\in{\cal V}}\inf\limits_{u\in{\cal
U}}H(t,x,p,u,v)=\inf\limits_{u\in{\cal U}}\sup\limits_{v\in{\cal
V}}H(t,x,p,u,v).$$

\end{lemma}

%                            lemma monotonie
%%%%%%%%%%%%%%%%%%%%%%%%%%%%%%%%%%%%%%%%%%%%%%%%%%%%%%%%%%%
\begin{lemma}\label{monotonie1}
The Hamiltonian $H$ satisfies  \bf{(H.2)} and \bf{(H.3)}.
\end{lemma}

The proof of the following results are similar to the proof of lemma
(\ref{monotonie}).
\begin{proposition}
1) For all $(u,v)\in {\cal U}\times {\cal V}$, let
$(Y^{u,v},Z^{u,v})$ be the solution of the BSDE with the generator
$(H(t,x,p,u_t,v_t),g_1(T,x_T)$ then $J(u,v)=Y_0^{u,v}.$\\
 2) Similarly, let $(Y^*,Z^*)$ be the solution of BSDE with
 generator $(H(t,x,p,u^*(t,x,p),v^*(t,x,p)),g_1(T,x_T))$ and define
 $(\widetilde{u},\widetilde{v})\in {\cal U}\times {\cal V}$ by
 $(\widetilde{u},\widetilde{v})=
 (u^*(t,x,Z^*_t),v^*(t,x,Z^*_t))_{t\leq T}$, then
 $J(\widetilde{u},\widetilde{v})=Y^*_0.$
\end{proposition}
\begin{theorem}
The strategy $(\widetilde{u},\widetilde{v})$ is a saddle-point for
the game.
\end{theorem}

\vskip 0.2cm

%%%%%%%%%%%%%%%%%%%%%%%%%%%%%%%%%%%%%%%%%%%%%%%%%%%%%%%%%%%
%                            proof lemma monotonie
%%%%%%%%%%%%%%%%%%%%%%%%%%%%%%%%%%%%%%%%%%%%%%%%%%%%%%%%%%%

%\begin{remark}In the same way we can deal with the stochastic zero-sum differential games and show that it has a value and a saddle point
%(one can see \cite{HLP95, HLP295} for more details).
%\end{remark}
%%%%%%%%%%%%%%%%%%%%%%%%%%%%%%%%%%%%%%%%%%%%%%%%%%%%%%%%%%%%%%%%%%%%%%%%%%%%%%%
%          References
%%%%%%%%%%%%%%%%%%%%%%%%%%%%%%%%%%%%%%%%%%%%%%%%%%%%%%%%%%%%%%%%%%%%%%%%%%%%%%%%

\end{document}